\documentclass[12pt]{amsart}
\usepackage{amsmath}
\usepackage{amsxtra}
\usepackage{amscd}
\usepackage{amsthm}
\usepackage{hyperref}
\usepackage[T1]{fontenc}
\usepackage{amsfonts}
\usepackage{amssymb}
\usepackage{eucal}
\usepackage[all]{xypic}
\usepackage{color} 
\usepackage[foot]{amsaddr}

\theoremstyle{theorem}
\newtheorem{thm}[subsection]{Theorem}

\newtheorem{lem}[subsection]{Lemma}
\newtheorem{prop}[subsection]{Proposition}

\theoremstyle{definition}

\newtheorem{defn}[subsection]{Definition}

\newtheorem{rem}[subsection]{Remark}

\newtheorem{example}[subsection]{Example}

\usepackage{yhmath}
\DeclareSymbolFont{largesymbols}{OMX}{yhex}{m}{n}
\DeclareMathAccent{\widetilde}{\mathord}{largesymbols}{"65}

\newcommand{\thmref}[1]{Theorem~\ref{#1}}

\newcommand{\lemref}[1]{Lemma~\ref{#1}}

\newcommand{\propref}[1]{Proposition~\ref{#1}}

\newcommand{\remref}[1]{Remark~\ref{#1}}
 
\newcommand{\exmref}[1]{Example~\ref{#1}}

\newcommand{\nc}{\newcommand}
\nc{\renc}{\renewcommand}
\nc{\ssec}{\subsection}
\nc{\sssec}{\subsubsection} 
\nc\ol{\overline}
\nc\wt{\widetilde}
\nc\wh{\widehat}
\nc\tboxtimes{\wt{\boxtimes}}

\renc{\d}{{\delta}}
\nc{\Aa}{{\mathbb{A}}}
\nc{\Bb}{{\mathbb{B}}}
 \nc{\Gg}{{\mathbb{G}}} 
\nc{\Hh}{{\mathbb{H}}}
 \nc{\Nn}{{\mathbb{N}}}
\nc{\Pp}{{\mathbb{P}}}
\nc{\Rr}{{\mathbb{R}}}
\nc{\BV}{{\mathbb{V}}}
\nc{\BW}{{\mathbb{W}}}
\nc{\Zz}{{\mathbb{Z}}}
\nc{\Qq}{{\mathbb{Q}}}
\nc{\Ss}{{\mathbb{S}}}
\nc{\Cc}{{\mathbb{C}}}
\nc{\Ff}{{\mathbb{F}}}

\nc{\CA}{{\mathcal{A}}}
\nc{\CB}{{\mathcal{B}}}

\nc{\CE}{{\mathcal{E}}}
\nc{\CF}{{\mathcal{F}}}

\nc{\CG}{{\mathcal{G}}}
\nc{\CL}{{\mathcal{L}}}
\nc{\CC}{{\mathcal{C}}}
\nc{\CM}{{\mathcal{M}}}

\nc{\CN}{{\mathcal{N}}}
\nc{\Oog}{{\mathbb{O}}}
\nc{\Oo}{{\mathcal{O}}}
\nc{\CP}{{\mathcal{P}}}
\nc{\CQ}{{\mathcal{Q}}}
\nc{\CR}{{\mathcal{R}}}
\nc{\CS}{{\mathcal{S}}}
\nc{\CT}{{\mathcal{T}}}
\nc{\CU}{{\mathcal{U}}}
\nc{\CV}{{\mathcal{V}}}
\nc{\CK}{{\mathcal{K}}}
\nc{\CW}{{\mathcal{W}}}
\nc{\CZ}{{\mathcal{Z}}}

\nc{\cM}{{\check{\mathcal M}}{}}
\nc{\csM}{{\check{\mathcal A}}{}}
\nc{\oM}{{\overset{\circ}{\mathcal M}}{}}
\nc{\obM}{{\overset{\circ}{\mathbf M}}{}}
\nc{\oCA}{{\overset{\circ}{\mathcal A}}{}}
\nc{\obA}{{\overset{\circ}{\mathbf A}}{}}
\nc{\ooM}{{\overset{\circ}{M}}{}}
\nc{\osM}{{\overset{\circ}{\mathsf M}}{}}
\nc{\vM}{{\overset{\bullet}{\mathcal M}}{}}
\nc{\nM}{{\underset{\bullet}{\mathcal M}}{}}
\nc{\oD}{{\overset{\circ}{\mathcal D}}{}}
\nc{\obD}{{\overset{\circ}{\mathbf D}}{}}
\nc{\oA}{{\overset{\circ}{\mathbb A}}{}}
\nc{\op}{{\overset{\bullet}{\mathbf p}}{}}
\nc{\cp}{{\overset{\circ}{\mathbf p}}{}}
\nc{\oU}{{\overset{\bullet}{\mathcal U}}{}}
\nc{\oZ}{{\overset{\circ}{\mathcal Z}}{}}
\nc{\ofZ}{{\overset{\circ}{\mathfrak Z}}{}}
\nc{\oF}{{\overset{\circ}{\fF}}}

\nc{\fa}{{\mathfrak{a}}}
\nc{\fb}{{\mathfrak{b}}}
\nc{\fg}{{\mathfrak{g}}}
\nc{\fgl}{{\mathfrak{gl}}}
\nc{\fh}{{\mathfrak{h}}}
\nc{\fj}{{\mathfrak{j}}}
\nc{\fm}{{\mathfrak{m}}}
\nc{\fn}{{\mathfrak{n}}}
\nc{\fu}{{\mathfrak{u}}}
\nc{\fp}{{\mathfrak{p}}}
\nc{\fr}{{\mathfrak{r}}}
\nc{\fs}{{\mathfrak{s}}}
\nc{\fsl}{{\mathfrak{sl}}}
\nc{\hsl}{{\widehat{\mathfrak{sl}}}}
\nc{\hgl}{{\widehat{\mathfrak{gl}}}}
\nc{\hg}{{\widehat{\mathfrak{g}}}}
\nc{\chg}{{\widehat{\mathfrak{g}}}{}^\vee}
\nc{\hn}{{\widehat{\mathfrak{n}}}}
\nc{\chn}{{\widehat{\mathfrak{n}}}{}^\vee}

\nc{\fA}{{\mathfrak{A}}}
\nc{\fB}{{\mathfrak{B}}}
\nc{\fD}{{\mathfrak{D}}}
\nc{\fE}{{\mathfrak{E}}}
\nc{\fF}{{\mathfrak{F}}}
\nc{\fG}{{\mathfrak{G}}}
\nc{\fK}{{\mathfrak{K}}}
\nc{\fL}{{\mathfrak{L}}}
\nc{\fM}{{\mathfrak{M}}}
\nc{\fN}{{\mathfrak{N}}}
\nc{\fP}{{\mathfrak{P}}}
\nc{\fU}{{\mathfrak{U}}}
\nc{\fV}{{\mathfrak{V}}}
\nc{\fZ}{{\mathfrak{Z}}}

\nc{\bb}{{\mathbf{b}}}
\nc{\bc}{{\mathbf{c}}}
\nc{\bd}{\partial}
\nc{\be}{{\mathbf{e}}}
\nc{\bj}{{\mathbf{j}}}
\nc{\bn}{{\mathbf{n}}}
\nc{\bp}{{\mathbf{p}}}
\nc{\bq}{{\mathbf{q}}}
\nc{\bF}{{\mathbf{F}}}
\nc{\bu}{{\mathbf{u}}}
\nc{\bv}{{\mathbf{v}}}
\nc{\bx}{{\mathbf{x}}}
\nc{\bs}{{\mathbf{s}}}
\nc{\by}{{\mathbf{y}}}
\nc{\bw}{{\mathbf{w}}}
\nc{\bA}{{\mathbf{A}}}
\nc{\bK}{{\mathbf{K}}}
\nc{\bI}{{\mathbf{I}}}
\nc{\bB}{{\mathbf{B}}}
\nc{\bG}{{\mathbf{G}}}
\nc{\bC}{{\mathbf{C}}}
\nc{\bD}{{\mathbf{D}}}
\nc{\bP}{{\mathbf{P}}}
\nc{\bH}{{\mathbf{H}}}
\nc{\bM}{{\mathbf{M}}}
\nc{\bN}{{\mathbf{N}}}
\nc{\bV}{{\mathbf{V}}}
\nc{\bU}{{\mathbf{U}}}
\nc{\bL}{{\mathbf{L}}}
\nc{\bT}{{\mathbf{T}}}
\nc{\bW}{{\mathbf{W}}}
\nc{\bX}{{\mathbf{X}}}
\nc{\bY}{{\mathbf{Y}}}
\nc{\bZ}{{\mathbf{Z}}}
\nc{\bS}{{\mathbf{S}}}

\nc{\sA}{{\mathsf{A}}}
\nc{\sB}{{\mathsf{B}}}
\nc{\sC}{{\mathsf{C}}}
\nc{\sD}{{\mathsf{D}}}
\nc{\sF}{{\mathsf{F}}}
\nc{\sG}{{\mathsf{G}}}
\nc{\sK}{{\mathsf{K}}}
\nc{\sM}{{\mathsf{M}}}
\nc{\sO}{{\mathsf{O}}}
\nc{\sQ}{{\mathsf{Q}}}
\nc{\sP}{{\mathsf{P}}}
\nc{\sZ}{{\mathsf{Z}}}
\nc{\sfp}{{\mathsf{p}}}
\nc{\sr}{{\mathsf{r}}}
\nc{\sg}{{\mathsf{g}}}
\nc{\sff}{{\mathsf{f}}}
\nc{\sfb}{{\mathsf{b}}}
\nc{\sfc}{{\mathsf{c}}}
\nc{\sd}{{\ltimes}}

\nc{\tH}{{\widetilde{H}}}
\nc{\tA}{{\widetilde{\mathbf{A}}}}
\nc{\tB}{{\widetilde{\mathcal{B}}}}
\nc{\tg}{{\widetilde{\mathfrak{g}}}}
\nc{\tG}{{\widetilde{G}}}
\nc{\TM}{{\widetilde{\mathbb{M}}}{}}
\nc{\tO}{{\widetilde{\mathsf{O}}}{}}
\nc{\tU}{\widetilde{U}}
\nc{\TZ}{{\tilde{Z}}}
\nc{\tx}{{\tilde{x}}}
\nc{\tq}{{\tilde{q}}}

\nc{\tfP}{{\widetilde{\mathfrak{P}}}{}}
\nc{\tz}{{\tilde{\zeta}}}
\nc{\tmu}{{\tilde{\mu}}}

 \nc{\vol}{{\mathop{\operatorname{\rm vol\,}}}}
 
 \nc{\colin}{{\mathop{\operatorname{\rm Colin\,}}}}
 \nc{\pic}{{\mathop{\operatorname{\rm Pic\,}}}}
 \nc{\disc}{{\mathop{\operatorname{\rm disc}}}}
 \nc{\Sym}{{\mathop{\operatorname{\rm Sym}}}}
 \nc{\Aut}{{\mathop{\operatorname{\rm Aut}}}}
 \nc{\Spec}{{\mathop{\operatorname{\rm Spec}}}}
 \nc{\spec}{{\mathop{\operatorname{\rm Spec}}}}
\nc{\Ker}{{\mathop{\operatorname{\rm Ker}}}}
 \nc{\dom}{{\mathop{\operatorname{\rm dom}}}}
\nc{\End}{{\mathop{\operatorname{\rm End}}}}
 \nc{\Hom}{\operatorname{\Hom}}
 \nc{\GL}{{\mathop{\operatorname{\rm GL}}}}
 \nc{\Id}{{\mathop{\operatorname{\rm Id}}}}
 \nc{\rk}{{\mathop{\operatorname{\rm rk}}}}
 \nc{\length}{{\mathop{\operatorname{\rm length}}}}
\nc{\supp}{{\mathop{\operatorname{\rm supp} \, }}}
\nc{\val}{{\rm val}}
 \nc{\res}{{\mathop{\operatorname{\rm res}}}}

\def\Ind#1#2#3{{#1} {\downarrow}_{#3} {#2} }
 
\def\tensor{{\otimes}}
\def\meet{\cap}
\def\union{\cup}

\def\si{\sigma}

 \def\m{\setminus}

\def\m{\smallsetminus}

\nc{\seq}[1]{\stackrel{#1}{\sim}}

\def\inv{^{-1}}

\def\beq#1{\begin{equation} \label{#1}}
\def\eeq{\end{equation}}

\def\prf{\begin{proof}}
\def\pv{\end{proof} }
 \def\eprf{\end{proof} }

 \def\lbl#1{ \label{#1} }

\def\kk{{\rm k}}
 \renc{\b}{{\beta}}

\def\Ind#1#2{#1\setbox0=\hbox{$#1x$}\kern\wd0\hbox to 0pt{\hss$#1\mid$\hss}
\lower.9\ht0\hbox to 0pt{\hss$#1\smile$\hss}\kern\wd0}

 \def\mS{\mathcal{S}}

\def\cti{\widetilde{\kk [t]}}
 \def\mI{\mathcal{I}}
 
 \renewcommand{\pic}{\operatorname{Pic}}

 \newcommand{\Th}{\operatorname{Th}}
 \renewcommand{\Aut}{\operatorname{Aut}}
 \renewcommand{\int}{int}

 \newcommand{\Div}{\operatorname{Div}}
 \renewcommand{\End}{\operatorname{End}}
 \newcommand{\Gal}{\operatorname{Gal}}
 \renewcommand{\supp}{\operatorname{supp}}
 \newcommand{\im}{\operatorname{im}}
 \newcommand{\tr}{\operatorname{tr}}
 
\begin{document}

\author{Taylor Dupuy$^1$} 
\address{$^1$Department of Mathematics and Statistics, University of Vermont, 82 University Place, Innovation Hall E220, Burlington VT 05405, United States}
\email{taylor.dupuy@gmail.com ($^1$corresponding author)}
\author{Ehud Hrushovski$^2$}
\address{$^2$Mathematical Institute, University of Oxford, Andrew Wiles Building, Radcliffe Observatory Quarter, Woodstock Road, Oxford
	OX2 6GG, United Kingdom}

\title{The theory of the entire algebraic functions}
\maketitle

\begin{abstract}
Let $A$ be the integral closure of the ring of polynomials $\Cc[t]$, within the field of algebraic functions in one variable. We show that $A$ interprets the ring of integers.
 This contrasts with the analogue for finite fields, proved to have a decidable theory in \cite{prestel-schmid}, 
 \cite{ vddries-macintyre}.
\end{abstract}

\begin{section}{Introduction}

For an integral domain $A$, let $A^{int}$ denote the integral closure of $A$ within the algebraic closure
of the field of fractions of $A$. 
Thus $\Zz^{int}$ is the ring of algebraic integers. 
Our motivation is 
the theorem of \cite{vddries}, that the theory of $\Zz^{int}$ is decidable (see also \cite{vddries-macintyre},
\cite{prestel-schmid}.)
There is a well-established analogy, going back to Leibniz and Euler,
between integers and polynomials; rational numbers and rational functions; algebraic numbers and algebraic functions.
The geometric counterpart of the arithmetic ring $\Zz^{int}$ is $\cti := \Cc[t]^{int}$ the ring of entire algebraic functions.\footnote{
An algebraic function is a multi-valued function $y$ of $t$, defined by a polynomial equation in $t$ and $y$. 
It is natural to call it \emph{entire} if $y$ never takes the value $\infty$ as long as $t \in \Cc$. 
This is equivalent to $y$ being integral over $\Cc[t]$, and for our purposes we simply define the ring of entire algebraic functions to be $\Cc[t]^{int}$.}
It is thus natural to ask as to the decidability of this ring. The answer turns out to go the other way. 
Recall that the {\em Kronecker dimension} of a field $k$ is the transcendence degree over the prime field, plus one if the characteristic is $0$. 

\begin{thm} \label{main} Let $\kk$ be a field of positive Kronecker dimension.
Let $\cti$ be the integral closure of $\kk [t]$ in ${\kk}(t)^{alg}$. 
Then the first-order theory $\Th(\cti,+,\cdot,t)$
interprets $\Th(\Nn,+,\cdot)$. 
\end{thm}

We will also show that $\Th(\cti,+,\cdot)$ is undecidable; indeed complete arithmetic is Turing-reducible to this theory, though we do not know if $\Qq$ is interpretable without the parameter $t$.

The only fields $\kk$ not in the scope of \thmref{main} are those algebraic over a finite field. For such fields,
it was indeed proved in \cite{vddries-macintyre} that $\kk[t]^{int}$ is is decidable.
 In Kronecker dimension zero, the fundamental principle is that of glueing. Given finitely many primes $p_1,\ldots,p_k$ of $\kk[t]$, and 
elements $a_1,\ldots,a_k$ of $\kk[t]$, in some 
 finite extension $L$ it is possible to find an integer 
that is close to $a_i$ at primes above $p_i$, and a unit at any other prime. This is no longer possible when $\kk$ is a non-algebraic field. When $C$ is an affine curve over
$\kk$, the obstruction to glueing on $\kk(C)$ in this sense is measured by the Picard group $\pic(C)$, namely the group of divisors modulo principal divisors; taken modulo torsion if one allows solutions in finite covers. In case $\bar{C} \m C$ has a single point, $\pic(C)$ can be identified with the Jacobian of the complete nonsingular curve $\bar{C}$;
in general it is a quotient of $J(\bar{C})$ by a countable subgroup (generated by the differences of points at infinity) and thus is not definable in the theory of algebraically closed fields.
At the limit over all curves $C$ one obtains in any case a $\Qq$-vector space $H$. We will see that $H$ is interpretable in $\cti$, and 
exhibit additional definable structure on this vector space; first a family of subspaces $H_b$ consisting of classes represented by divisors supported above a finite subset of some curve; then finer structure leading to definability of the geometry of linear dependence on $H$. 
It remains a very interesting open problem 
to find a decidable setting capturing the essential geometry here. 

With regard to the specific structure $\cti$ in any case, the coefficient fields $\Qq$ associated with Picard groups of elliptic curves become interpretable, implying a strong undecidability, and in this way joining a large body of literature starting with \cite{raphael-robinson}. Unlike our study, the emphasis is usually on subrings of number fields or finitely generated function fields. 
Works including \cite{denef}, \cite{pheidas}, \cite{poonen}, \cite{MR}, \cite{LMB}, \cite{LMBsh}, \cite{jk} make use in this context of elliptic multiplication, that also plays a large if somewhat different role here. 
\end{section}

\section{Affine curves, principal divisors and supports}

Let ${\kk}$ be an algebraically closed field. We consider nonsingular affine curves $C$ over $\kk$.
(They will arise as the spectrum of $\cti \meet L$ for various finite extension fields $L$ of $\kk(t)$.) 
By the {\em genus} $g$ of an affine curve $C$, we mean here the genus of the associated smooth projective curve $\bar{C}$.
We have a natural map $\si: C^g \to \pic(C)$, where $\pic(C)$ is the Picard group of $C$.
 
Let $C$ be an affine curve.
Let $\Bb(C)$ be the Boolean algebra of finite or cofinite subsets of $C({\kk})$, ${{\Bb_f}}(C)$ the ideal of finite subsets,
Let $\Div(C)$ be the free Abelian group on $C({\kk})$, viewed as the group of functions $C({\kk}) \to \Zz$ with finite support.
For $f \in \kk(C)^*$, define the affine divisor $(f)_{C}$ to be the element of $\Div(C)$ with $(f)_{C}(p) : = v_p(f)$, $v_p$ being the valuation corresponding to $p$. 
Let $D_{pr}(C)$ be the group of affine divisors of elements of
$ \kk(C)^*$.
We define $\pic(C)$ by the exact sequence:
\[ 0 \to D_{pr}(C) \to \Div(C) \to \pic(C) \to 0,\] 
 so $\pic(C)$, the Picard group of $C$, is the quotient of $\Div(C)$ by the principal divisors.
 Let $\rho=\rho_C: \Div(C) \to \pic(C)$ be the natural map.
 
\ssec{Warmup: support within a single curve} 
 
For $b \in \Bb_f(C)$, let $\Div(C)_b$ be the free Abelian group on $b$.
An element of $\Div(C)_b$ will be said to be {\em supported} on $b$. 
An element $a$ of $\pic(C)$ is said to be supported on a subset $\{a_1,\ldots,a_n\}$ of the curve $C$ if some multiple $na$ is represented by an expression $\sum m_i a_i$ (modulo principal divisors.) 
 
In $\cti$ we do not have direct access to a single curve, and cannot discuss this notion of support (except at the limit). The lemma below will thus not be used, and is meant merely as an illustration. It shows {\em at the level of a single curve} how to define $\Qq$-linear dependence
 on the Picard group, given the notion of support on a fixed curve. We assume for simplicity that $\kk$ has infinite transcendence degree over the prime field.

\begin{example}\label{warmup} 
 Assume $\kk$ has infinite transcendence degree over the prime field. 
Let $C$ be a smooth affine curve over $\kk$, $a,a' \in \pic(C)$. Assume $a'$ is supported on $b$ whenever $a$ is supported on $b$. Then $a' \in \Qq a \in \Qq\otimes \pic(C)$. 
 
 \end{example}
\prf We have to show that $a' \in \Qq a$. Let $k_0$ be a finite transcendence degree subfield of $\kk$, with $C$ defined over $k_0$. Let $g$ be
the genus of the complete nonsingular curve $\bar{C}$ associated to $C$; so $\pic(C)$ is a quotient group of the Jacboian $J$ of $\bar{C}(\kk)$, with finitely generated kernel.
 Pick generic elements $b_1,\ldots,b_g$ of $C$, over $k_0(a,a')$.
Let $c_1,\ldots,c_g$ be elements of $C$ such that $a+\sum_{i=1}^g b_i= \sum_{i=1}^g c_i$ as elements of $\pic(C)$. Then in $\Qq \tensor \pic(C)$, (since the same is true of $a$,) also $a' $ is a $\Qq$-linear combination of $a,b_1,\ldots,a_g$ as well as of $c_1,\ldots,c_g$, say $a' = \sum _{i=1}^g \gamma_i c_i$ 
 If $a' =0$ there is nothing to prove; otherwise, some $\gamma_i \neq 0$, say for $i=1$; replacing $a'$ by $a'/\gamma_1$ we may assume $\gamma_1 = 1$.
 Thus $a-a' + \sum_{i=1}^g b_i = \sum _{i=2}^g (1-\gamma_i) c_i$ in $\Qq \tensor \pic(C)$. 
Thus in $J$, 
 $ \sum_{i=1}^g b_i$ lies in a translate (defined over $k_0(a,a')^{alg}$) of $P$, with $P$ the $g-1$-fold sums of elements of the curve in $\bar{C}$.
This contradicts the genericity of $b$. 
 \eprf

 The challenge we will face beyond this is that even a single element $c$ of $C$, considered as a support, lifts to a bigger finite subset $c_1',\ldots,c_k'$ of a covering curve $C'$; 
 to say that an element is supported on $c$ includes the possibility that it is supported on some subset of $\{c_1',\ldots,c_k'\}$; and we cannot choose parameters generic to all
 possible $C'$, but rather must consider ramified coverings that depend on unknown parameters. The phenomenon seen in \exmref{warmup} remains true all the same; but to prove it we will have to take a closer look using trace maps for Galois group actions. We prepare this with a lemma.

 \begin{lem} \lbl{ic6a} Let $F$ be a finitely generated field of positive Kronecker dimension. 
 Let $C,D$ be smooth affine curves over 
 $F$, and $p: D \to C$ a finite morphism. 
 Let $a \in \pic(C)$, and let $d$ be a tuple from $F^{alg}$.
 Then there exist 
 $b=(b_1,\ldots,b_g) $ and $c=(c_1,\ldots,c_g)$ in $ C(F^{alg})^g$ such that $\si(b)+ a = \si(c)$, and 
 $\Aut(F^{alg} / F(a,b,c,d))$ acts transitively on $\Pi_{i=1}^g p \inv(b_i)$, as well as on 
 $\Pi_{i=1}^g p \inv (c_i)$. \end{lem}

 \prf Let $p_1: C \to \Pp^1$ be a dominant morphism. Let $F'=F(b_1',\ldots,b_g')$ with
 $b'=(b_1',\ldots,b_g')$ a generic element of $C^g$. Let $t'_i=p_1(b_i')$; then 
 $F(t'_1,\ldots,t'_g)$ is a purely transcendental extension of $F$. Let $c=(c_1',\ldots,c_g') \in C(F')^g$
 be such that $\si(b') + a = \si(c')$. 

 If $e_i', e_i'' \in p \inv(b_i)$, then $(e_1',\ldots,e_g')$
 is a generic element of $D^g$ over $F(d)$ as is $(e_1'',\ldots,e_g'')$; so there exists a field automorphism
 fixing $F(d)$ and 
 taking each $e_i' \mapsto e_i''$; this field automorphism fixes $F(b',c',d)$. Similarly,
 $Aut(F(b',c',d)^{alg} / F(b',c',d))$ acts transitively on $\Pi_{i=1}^g p \inv (c'_i)$.

 Let $L'$ be a finite normal field extension of $F(t'_1,\ldots,t'_g,d)$ and such that: 
 \begin{enumerate}
 \item Each $b_i'$ and $c_i'$ lie in $C(L')$; we have $\si(b')+a=\si(c')$, and $p_1(b_i')=t_i$; moreover, 
 \item all elements of $D(F(t'_1,\ldots,t'_g)^{alg})$ lying above some $b_i'$ or $c_i'$, belong to ${L'}$. 
 \item $\Aut({L'}/F(b',c',d))$ acts transitively on $\Pi_{i=1}^g p \inv (b'_i)$ and on $\Pi_{i=1}^g p \inv (c'_i)$.
 \end{enumerate}
 
 By \cite{fj} Theorem 13.4.2,
 $F$ is a Hilbertian field, and in fact we may embed $F(t'_1,\ldots,t'_g)$ in an elementary extension $F^*$ of $F$,
 in such a way that $F^*$ is a regular field extension of $F(t'_1,\ldots,t'_g)$,
 i.e. it is linearly disjoint from $F(t'_1,\ldots,t'_g)^{alg}$ over $F(t'_1,\ldots,t'_g)$. 
 It follows
 that the compositum ${L'}F^*$ is a normal extension of $F^*$ with automorphism group $G=\Aut(L'F^*/F^*)=\Aut(L'/F(t_1',\ldots,t'_g))$. 
 
 Recall the set of field extensions $E$ of a given degree of any field $F$ is interpretable in $F$, along with the automorphism group $\Aut(E/F)$. 
Applying this to $L'F^* / F^*$, it is easy to write a formula $\psi(u_1,\ldots,u_g)$ in the language of rings,
 such that $F^* \models \psi(t_1',\ldots,t'_g)$, and such that $\psi$ asserts the 
 existence of a finite Galois field extension having the properties above.
 Since $F \prec F^*$, there exist $t_1,\ldots,t_g \in F$ with $F \models \psi(t_1,\ldots,t_g)$; hence
 there exist a finite normal field extension $L$ and elements $b_i,c_i \in C(L)$ having the properties (1-3)
 above. 
 \eprf 

\section{Definable systems of class groups} 

Let ${\kk}$ be an algebraically closed field. 
Fix an algebraic closure ${\kk}(t)^{alg}$ of the field of rational functions ${\kk}(t)$. 
Let $\cti$ be the integral closure of ${\kk}[t]$ in ${\kk}(t)^{alg}$. 
Let $\mI$ be the set of finite extension fields of ${\kk}(t)$, contained in ${\kk}(t)^{alg}$.
This is a directed partially ordered set, under inclusion. We will also denote the field $i$ by $k_i$. 
Each $i \in \mI$ is the function field of a smooth projective curve $\bar{C}_i$, over ${\kk}$. 
Let $C_i$ be the Zariski open subcurve of $\bar{C_i}$, whose affine coordiante ring is $\kk_i \meet \cti$. 
 
Given $i \leq j \in \mI$ there is a unique morphism $p_{ji}: C_j \to C_i$ of curves, corresponding to the given inclusion of function fields; the $p_{ji}$ form a directed system.

Consider now a pair $i \leq j \in \mI$.
Define a homomorphism of partially ordered groups $ {p^{ij}}: \Div(C_i) \to \Div(C_j)$ by the formula 
\beq{dci}{ {p^{ij}}(d)(q) = r(q) d(p_{ji}(q)),} \eeq
where $r(q)$ is the ramification degree of $j$ over $i$ with respect to the valuation $q$.
The multiplicities $r(q)$ are included so as to have, for any $f \in {\kk}(C_i)$, 
\[p^{ij} (f)_{C_j} = (f \circ p_{ji})_{C_i},\] 
so $ p^{ij}(D_{pr}(C_i)) \subseteq D_{pr}(C_j)$.
Thus ${p^{ij}}$ also induces a homomorphism 
\beq{jfin} { p^{ij}: \pic(C_i) \to \pic(C_j). } \eeq

Finally, we have a natural embedding
\beq{jfin3} p^{ij}: {{\Bb_f}}(C_i) \to {{\Bb_f}}(C_j),  \eeq
 namely $p^{ij}(y) = p_{ji}^{-1} (y)$. It is compatible with $p^{ij}: \Div(C_i) \to \Div(C_j)$
and the support maps $\supp: \Div(C_i) \to {{\Bb_f}}(C_i)$, taking an element to its support. 

Let $\Bb,{{\Bb_f}}, D, {{H}}, H_b$ be the direct limits along $\mI$ of $\Bb(C_i), {{\Bb_f}}(C_i) ,\Div(C_i), \pic(C_i), H_b(C_i)$. The transition maps from $i$ to $j$
are \eqref{dci},\eqref{jfin}, \eqref{jfin3}. 
$\Bb$ is a Boolean algebra with maximal ideal ${{\Bb_f}}$.
 These direct limits come with the maps $p^{i \infty}: \Div(C_i) \to D$, 
$p^{i \infty}: H_b(C_i) \to H_b$, etc.

\begin{lem} \lbl{ic2} $D,H$ and the $H_b$ are torsion-free divisible abelian groups. \end{lem} 

\prf Each $\Div(C)$ is torsion-free, hence so is $D$. And
 any element of $D$ attains an $n$'th root with the same support in some (sufficiently ramified) covering.
 (If as in \lemref{ic1} we write $d= \min( (f),(g))$, then $\frac{1}{n} d = \min (f^{1/n},g^{1/n})$, so it suffices to go to a cover of $C$
 on which $f^{1/n},g^{1/n}$ are regular functions.)
It follows that $H$ is divisible.

Let us show that $H$ is torsion-free. We have to show that for any $i$, with $C=C_i$, any torsion element of $\pic(C)$ maps to $0$ in $\pic(C')$ for an appropriate covering curve $C' \to C$. Let $\tau \in \Div(C)$ represent a torsion element of $\pic(C)$; 
so there is a function $f \in {\kk}(C)$ with $(f)_C=m \tau$. Let $k_j$ contain $k_i(g)$ with $g^m=f$. 
Then $(g^m) = m p^{ij}(\tau)$ so $(g)= p^{ij}(t)$, thus $\tau$ pulls back to $0$ in $\pic(C')$. 
A similar argument works for $H/H_b$, for any $b \in {\Bb_f}$: if $\tau$ represents an $m$-torsion element of $H/H_b$, 
we have $(f)= m \tau + \upsilon$ with $\upsilon$ supported above $b$, and so 
$(g) = \tau + \upsilon'$ with $\upsilon'$ still supported above $b$. 
 
Since $H/H_b$ is torsion-free, it follows that $H_b$ is divisible too.
\eprf

Keeping in mind the $\Qq$-vector spaces at the limit, our interest in the approximations $\pic(C)$ will always be modulo torsion.
For a subgroup $A$ of $\pic(C)$ we will write $A \approx (0)$ if $A$ is a torsion group.

We fix an element of $ \mI$ corresponding to an affine curve $C_1$, and denote it by $1$; the corresponding function field is $k_1$. 
For any $ j \geq 1$ we also have a morphism $p_{1j*}: \Div(C_j) \to \Div(C_1)$ in the opposite direction, where for $e: C_j(k) \to \Zz$ the map $p_{1j*}e:C_1(k) \to \Zz$ is defined by 
\[ p_{1j*}(e)(p) = \sum_{p_{j,1}(q)=p} e(q). \]
Note that the composition $p_{1j*} \circ {p^{1j}} $ is multiplication by $n=[k_j:k_1]$. 

For a principal divisor $e = (f)_{C_j}$, we have $p_{1j*}(e) = (Nf)_{C_1}$, where $Nf$ is the image of $f$ under the norm map from $\kk(C_j)$ to
$\kk(C_1)$. 
Since $p_{ij*}$ and $p^{1j}$ respect the groups of principal divisors, we have induced homomorphisms between 
 $\pic(C_j)$ and $\pic(C_1)$ denoted by the same letters.
We still have $p_{1j*} {p^{1j}} = [\cdot n]$ with $n=[k_j:k_1]$. 
Thus on $\pic(C_j)$ we have 
\[ \ker p_{1j*} \meet \im p^{1j} \approx (0). \]
On the other hand $\ker p_{1j*} + \im p^{1j} = \pic(C_j)$ (as in \lemref{ic3} below.) 
Note that for $j \leq k$, $p^{jk}(\ker p_{1j*}) \subset \ker p_{1k*}$.

Let ${H}_{1} = p^{1\infty}(\pic(C_1)) \cong \Qq \tensor \pic(C_1)$ be the image of $\pic(C_1)$
in ${{H}}$, and let ${{H_1}^{\perp}} = \varinjlim_j \ker p_{1j*}$ . Then ${{H_1}^{\perp}}, H_1$ are complementary $\Qq$-subspaces
of $ {H}$.

Fix $b \in {{\Bb_f}}$. An element of ${{H}}$ has many representatives in $D$. We let $H_b$ 
be the set of elements of ${{H}}$ having some representative supported on $b$. (Sometimes we will use the notation $H(b)$.)

\begin{defn} 
For $b \in {\Bb_f}(C_i)$, $i \leq j$, let ${H}(j;b) = \rho(D(j;b))$, where $D(j;b)$ is the set of elements of $\Div(C_j)$ supported on (a subset of) $p^{ij} b$. Let
${H}_b$ be the limit of the ${H}(j;b)$ over $j$. 
\end{defn}

The $H(j;b)$ are our fundamental geometric objects. 

\ssec{ The subspaces ${H}(b)$.} 
Say $b=\{\b_1,\ldots,\b_k\} \subset C_1({\kk})$. Let ${H}_1(b)$ be the subspace of ${H}$ generated by the images of $\b_1,\dots,\b_k$.
Let ${{H_1^{\perp}}}(b)$ be the subspace of ${H}$ generated by the differences $x-y$ where $p_{j1}(x)=p_{j1}(y)=\b_i$ for some $i \leq k$. 
Note that ${H}_1(b) \leq {H}_1$, while ${{H_1^{\perp}}}(b) \leq {{H_1}^{\perp}}$.

\begin{lem} \lbl{ic3} ${H}(b) = {H}_1(b) \oplus {{H_1^{\perp}}}(b)$ .
\end{lem}

\prf The sum is direct since ${H}_1 \meet {H_1}^{\perp} = 0$. To show that $ {H}_1(b) \oplus {{H_1^{\perp}}}(b) =H(b)$, consider an element of $\Div(C_j)$ supported above $b$. It is a linear combination of singleton elements
of $C_j$, supported above some $\beta=\beta_l$. So it suffices to show that
such an $e$ lies in $ {H}_1(b) + {{H_1^{\perp}}}(b) $. Say $n=[C_j:C_1]$. 
Then $ne = p^{1j}(\beta) - \sum_{p_{j1}(y)=\beta} (y-e)$, and each $y-e \in {{H_1^{\perp}}}(b)$. This
exhibits $ne$ as the difference of an element of ${H}_1(b)$ and one of $ {{H_1^{\perp}}}(b) $. Since
by definition these are $\Qq$-subspaces, we have also $e \in {H}_1(b) + {{H_1^{\perp}}}(b) $.
\eprf

 We now take a closer look at the interaction of the subspaces ${H}(b)$. For simplicity, we state the lemma first in genus one. (This is also the case that will be used for the interpretation of $\Qq$.)

\begin{lem} \lbl{ic6} Let $C_1$ be an affine curve of genus one over a field $\kk=\kk^{alg}$ of positive Kronecker dimension. 
Fix $j$; so we have a covering $p=p_j: C_j \to C_1$.
 
 Let $\b_1 \in C_1(k)$. Then there exists $\b_2 \in C_1(k ) \ $ such that,
 letting $\b_3=\b_1-\b_2$, statements 
 (1-3) below
 hold. 
 \begin{enumerate}
 \item ${H}(j;\beta_1) \meet ({{H_1^{\perp}}}(j;\beta_2) + {{H_1^{\perp}}}(j;\beta_3)) \approx (0)$ 
\item $ {{H_1^{\perp}}}(j;\beta_1) \meet ({H}(j;\beta_2) + {H}(j;\beta_3)) \approx (0)$ 
\item $ p^{j \infty} {{H_1^{\perp}}}(j;\beta_1) \meet ({H}( \beta_2) + {H}(\beta_3)) = (0) $
 \end{enumerate}
 \end{lem}

\prf We may choose a finitely generated subfield $k_0$ of $k$ such that $k_0$ has 
positive Kronecker dimension, $C_j,C_1,p$ are defined over $k_0$, and $\b_1 \in C_1(k_0)$.

 Let ${A}=\pic(C_1)$ and let ${B}=\pic(C_j)$ be the Picard group of $C_1,C_j$. 
 Choose $\beta_2,\beta_3$ as in \lemref{ic6a} (i.e. $\beta_2=b,\beta_3=c$ with respect to $ a=\beta_1, D=C_j, d$ an enumeration of 
 $p \inv(a)$, in the notation of that lemma; note here $g=1$).
 
 Let $a_1 \in B(k)$ represent an element of $ {H}(j;\beta_1)$; 
 in other words, $a_1$ is represented by a cycle on $C_j$ whose support lies above $\beta_1$. 
 Let $a_1=a_2+a_3$ such that for $i=2,3$, we have $a_i$ supported above $\beta_i$, 
 and $p(a_i)$ is supported above $\delta$. 
 We will show that $a_1$ is torsion.
 Note that $a_1 \in B(k_0^{alg})$, since 
 $ \beta_1 \subset C_j(k_0^{alg})$, and $a_1$ is supported on points above this finite set. 
 Let $K=k_0^{alg}(\b_2) = k_0^{alg}(\b_3)$, $L$ the Galois hull over $K$
 of $K(a_2)=K(a_3)$, and $G=\Aut(L/K)$, $m=|G|$. Then $G$ acts on $B(L)$ and on $A(L)$, and we have a trace
 map (sum of conjugates) $\tr: B(L) \to B(K)$. Since $G$ fixes $a_1$ we have $\tr(a_1)=ma_1$.
 On the other hand for $i=2,3$, any two elements of $p_{j1}\inv(\b_i)$ are $\Aut(L/K)$-conjugate; thus their difference has $G$-trace zero. Since $a_2$ and $a_3$ 
 (being in $H_1^\perp$) 
 are sums of such differences, 
 we have $\tr(a_2)=\tr(a_3)=0$. 
 Thus $ma_1=0$. This proves (1).
 
 We now deduce (2,3) from (1). 
 For (2), assume $a_1 \in {{H_1^{\perp}}}(j;\beta_1)$. We apply the operator $\operatorname{id} - p^{1j} p_{j1 *} $ where $\operatorname{id}$ denotes the identity;
 it leaves $a_1$ fixed since $a_1 \in {{H_1^{\perp}}}(j;\beta_1)$, so that $p_{j1*}(a_1)=0$; and for $i=2,3$ takes $a_i$
 to elements of ${{H_1^{\perp}}}(j;\beta_i)$, so that the previous paragraph applies. 
 
For the last point, (3), it suffice to show for any ${j'} \geq j$, that 
 \[ p^{j {j'}} {{H_1^{\perp}}}(j;\beta_1) \meet ({H}({j'}; \beta_2) + {H}({j'};\beta_3)) \approx (0). \]
 Let $a_1 \in {{H_1^{\perp}}}(j;\beta_1)$, $a_i \in H({j'};\beta_i)$ for $i =2,3$ and suppose
 $p^{j{j'}}(a_1) = a_2+a_3$. Applying $p_{{j'}j*}$, and using that $p_{{j'}j*}p^{j{j'}}(a_1)=da_1$ for appropriate $d$,
 we see that $d a_1 \in H(j;\beta_2)+H(j;\beta_3)$
 so $a_1$ is torsion. 

 \eprf

 Observe that \lemref{ic6} (3) 
 does not go as far as asserting that $ {{H_1^{\perp}}}(j';\beta_1) \meet ({H}(j';\beta_2) + {H}(j';\beta_3)) \approx (0)$.
 This is because even if $\beta_1$ is chosen to be generic over $k_0$, it cannot be chosen
 generic over all $j'$.

Here is the more general statement of \lemref{ic6}; the proof is identical.

\begin{lem} \label{ic6g} Let $C_1$ have genus $g \geq 1$. Let $p=p_j: C_j \to C_1$ be defined over a subfield $k_0$
of $k$. 
Let $\b_{1,1},\ldots,\b_{1,m} \in C_1(k_0)$, $\beta_1 =( \b_{1,1},\ldots,\b_{1,m} )$. Then there exist two $g$-tuples $\b_2 = \{\b_{2,1},\ldots,\b_{2,g}\}$ 
and 
$\beta_3=(\b_{3,1},\ldots,\b_{3,g})$ of elements of $C_1$ such that in the Picard group $\pic(C_1)$ we have
\[ \sum_{\nu=1}^m \b_{1,\nu} = \sum_{\mu=1}^g \b_{2,\mu} + \sum_{\mu=1}^g \b_{3,\mu}, \]
and \[ p^{j \infty} {{H_1^{\perp}}}(j;\beta_1) \meet ({H}( \beta_2) + {H}(\beta_3)) = (0).\]
\end{lem}

 \begin{section}{Interpretations}

We consider the structure $\mS$ consisting of the Boolean algebra $\Bb$, the ideal ${{\Bb_f}}$, 
the partially ordered group $D$, 
the support map 
$D \to {{\Bb_f}}$, the group ${{H}}$ and the homomorphism $\rho: D \to {{H}}$. 

We will show below that $\mS$ is interpretable in $\cti$. The Boolean ingredients appear already in \cite{vddries}, \cite{vddries-macintyre}, \cite{prestel-schmid}, 
 and the definability of finitely generated ideals and of their radicals is a key point there.

\begin{lem} \lbl{ic1} Any finitely generated ideal of $\cti$ is $2$-generated as an ideal. 
For any affine curve $C$, any function $C (\kk) \to \Nn$ with finite support has the form $p \mapsto \min(v_p(f),v_p(g))$
for some $f,g \in \cti$. 
\end{lem}

\prf This lemma is classical, see e.g. \cite{vddries-macintyre}; we give a proof in the present geometric setting.
 Any ideal $I$ of the Dedekind domain $R$ has the form $\{f: v_{p_k}(f) \geq m_k ,k=1,\ldots,m\}$ for some $p_1,\ldots,p_k \in C({\kk})$ and $m_1,\ldots,m_k \in \Nn$. Since any two of the ideals $p_i^{m_i}$ generate the ideal $R$, 
using the Chinese remainder theorem, we may find $g \in R$
with $v_{p_k}(g)=m_k$. Now $g$ may have additional zeroes at some additional points $q_1,\ldots,q_l$. But by the independence of valuations,
there exists $h$ with $v_{p_i}(h) \geq m_i$ and $v_{q_i}(h)=0$. Now it is clear that $g,h$ generate $I$. The second statement is proved in the same way. \eprf

\begin{lem} \lbl{ic2b} $\mS$ can be interpreted in the ring $\cti$. \end{lem}
 \prf Let $A=\cti$. Finitely generated ideals are uniformly definable as by \lemref{ic1} they are all $2$-generated. As for radicals: 
 since all prime ideals of $A$ are maximal, the radical of an ideal, intersection of all prime ideals containing it,
 coincides with the Jacobson radical; the standard formula for the Jacobson radical of a definable ideal $I$
 shows that $\sqrt I$ is definable too.

 ${{\Bb_f}}$ can be identified with the set of radicals of finitely generated 
ideals of $\cti$.
Intersection in ${{\Bb_f}}$ corresponds to the radical of the sum of ideals; union to the
intersection. The difference of two elements of ${{\Bb_f}}$ corresponds to the operation
of (radical) ideals, mapping $(I,J)$ to $ \{x: Jx \subset I \}$. Given ${{\Bb_f}}$
with this structure, the Boolean algebra $\Bb$ can be easily interpreted by adding a formal
element $1$ and defining the obvious structure on $ {{\Bb_f}} \union \{1-a: a \in {{\Bb_f}}\}$.

Let $D^+$ be the semigroup of non-negative elements of $D$. 
Define a map $\alpha: A^2 \to D^+$, $\alpha(a,b)(p)=\min(v_p(a),v_p(b))$; more precisely, define this at the level of each sufficiently large field $i$ of the limit system, and check compatibility. By \lemref{ic1},
$\alpha$ is surjective. Let us show that equality and the ordered semigroup operations 
$+,\min$ pull back to definable sets on $A^{4}$. Indeed $\min$ corresponds to 
the sum of ideals, i.e. $\alpha(a'',b'')=\min(\alpha(a,b),\alpha(a',b'))$ iff $A(a,b,a',b')=A(a'',b'')$.
Equality corresponds to equality of ideals. $+$ corresponds to product of ideals:
$\min(a,b)+\min(c,d) = \min(ac,ad,bc,bd)$.
 The lattice-ordered group
$D$ is now easily definable as the set of differences $a-b$ with $a,b \in D^+$.
 
 Observe that $D$ is interpreted along with the structural map $\alpha: A^2 \to D^+ \subseteq D$.
 In particular we can define $\beta(a)=\alpha(a,a)$ so that 
$\beta(a)(p)=v_p(a)$. We define ${{H}}$ as $D$ modulo the image of $\beta$.

The support map $D^+ \to {\Bb_f}$ is obtained by factoring through $D$ the map $\min: A^2 \to {\Bb_f}$:
\[ \min((a_1),(a_2)) \mapsto \sqrt{(a_1,a_2)}, \]
where $\sqrt{(a_1,a_2)} $ is the radical ideal generated by $a_1,a_2$. This then extends to $D$, the support of $d$
being the union of supports of $d^+= \min(d,0)$ and $d^-= -\max(d,0)$.

\eprf
 
\begin{rem} \label{connectiontok} Let $\kk_1$ be a subfield of $\kk$, let $C$ be a $\kk_1$-definable curve, namely the spectrum of
$\Oo:= \cti \meet L$, where $L$ is a finite extension field of $\kk(t)$. 
Let $c \in C(\kk_1)$ and let $d$ be the image of $c$ in $D$. Then $d$ is $\Oo$-definable in $\cti$. 
Indeed the point $c$ of $C$ corresponds to a homomorphism $\cti \to \kk_1$, whose kernel generates the radical ideal corresponding to $d$. Hence the image
of $c$ in $H$ (namely $\rho(d)$) is also $\Oo$-definable in $\cti$.
\end{rem}

\noindent {\em Assume $k$ has positive Kronecker dimension.} We proceed to interpret $\Qq$ in $\mS$. 
 
For a curve $C_1$ over $\kk$, recall $H_1$ is the image of $\pic(C_1)$ in $H$.
 
 We have a definable family $\Xi_1$ of $\Qq$-subspaces of ${H}$, namely all the ones of the form ${H}(b)$ for $b \in {\Bb_f}$. 

So $\Xi = \{ U+V: U,V \in \Xi_1 \}$ is also a uniformly definable family. Actually $\Xi=\Xi_1$ since
$H(b \union b')=H(b)+H(b')$.

For $a \in {H}$, let 
$\Phi(a)$ denote the intersection of all $U \in \Xi$ with $a \in U$. So $\Phi(a)$ is definable uniformly in $a$; it is a $\Qq$-space; and $a \in \Phi(a)$. 

\begin{lem} \label{1d} For any nonzero $a \in H$ represented by an irreducible divisor on $C_1$, $\Phi(a)$ is the one-dimensional subspace $\Qq a$ of ${H}$.
\end{lem}
 \prf
 Let $a \in H$ be represented by an irreducible divisor $a_0$ of $C_1(\kk)$. Let $b \in \Phi(a)$. We have to show that $b \in \Qq a$. 
 During the proof we may replace $b$ by a positive multiple when needed.
 
 Since $a$ is supported on $a_0$ and $b \in \Phi(a)$, we have $b \in {H}(\{a_0\})$. By \lemref{ic3}, $b \in {H}_1(\{a_0\}) \oplus {{H_1^{\perp}}}(\{a_0\})$
 Subtracting a scalar multiple of $a$ from $b$, we may assume $b \in {{H_1^{\perp}}}(\{a_0\})$. So for some $j$ we have
$b \in p^{j \infty} {{H_1^{\perp}}}(j;a_0)$.


Let $\beta_2,\beta_3$ be as in \lemref{ic6g}.
 Then $a \in {H}(\b_2)+{H}(\b_3)$. Since $b \in \Phi(a)$ we have (*) $b \in {H}(\b_2)+{H}(\b_3)$.
 But then by \lemref{ic6g} we have $b=0$, so $b \in \Qq a$. 
 
 \eprf 
 \begin{rem} \label{1dr} The same proof shows that if $a= \sum_{i=1}^n a_i$ in $H$, with $a_i$ represented by an irreducible divisor on a curve $C$,
 then $\Phi(a)$ is contained in the span of $a_1,\ldots, a_n$. Since any element is the image of some element of $\pic(C)$ for some curve $C$ in the directed system,
 it follows that for any $a \in H$, $\Phi(a)$ is finite-dimensional.
\end{rem}
 
 Note that if $C_1$ has genus $1$, then any element of $\pic(C_1)$ is represented by the irreducible divisor consisting of a single point of $C_1$.

\begin{proof}[Proof of \thmref{main}] 

Note that in \thmref{main} we may assume $\kk$ is algebraically closed, since changing $\kk$ to $\kk^{alg}$
does not change the ring $\kk[t]^{int}$.
We work in $\cti$ with parameter $t$. By \lemref{ic2b}, we have also the structure $\mathcal{S}$ at our disposition.

Fix a smooth affine curve $C_1$ of genus $1$ over the prime field $\kk_0$ defined by a polynomial $F(x,t)=0$ in full Weierstrass form, quadratic in the $x$ variable and hence Galois over the $t$-line.
Let $b_1,b_2$ be points of $C_1(\kk_1)$ with linearly independent images $a_1,a_2$ in $H$; where $\kk_1 \leq \kk$ is a finite extension of $\kk_0$.
Let $L = \kk_1(C_1)$, a finite extension of $\kk_1(t)$ corresponding to the above morphism from $C_1$ to the $t$-line. 
 We will at first use parameters not only from $\kk_0(t)$ but also from $L$. Below we describe how to avoid $L$-parameters and restrict to $\kk_0(t)$.

Note that $a_1,a_2$ are definable elements of $H$ in $\cti$, over parameters from $L$ (\remref{connectiontok}).
Let $ V=\Qq a_1 + \Qq a_2 $. Then $V=\Phi(a_1)+\Phi(a_2)$ is definable, since $+$ is definable on ${H}$.
Letting $b$ vary in $H$ and considering $\Phi(b) \meet V$, we obtain a family of subspaces of V, including all one-dimensional subspaces
of $V$. Since addition on $V$ is definable, we can uniformly define all affine subspaces of $V$ too. Hence colinearity is definable, as a ternary relation
$\colin$. Now it is easy to interpret the coefficient field $\Qq$ in $(V,+,0,a_1,a_2, \colin)$. 
This is one direction of the Cartesian algebra / geometry connection (later encapsulated in the fundamental theorem of projective geometry.)
We can take the universe of $\Qq$ to be the $x$-axis,
 i.e. the line through $(0,0)$ and $(1,0)$. The bijection with the $y$ axis can be defined by mapping $a$
to $b$ if the line through $(a,b)$ is parallel to the line through the basis elements $(1,0)$, $(0,1)$. 
Multiplication $a \cdot b$ can be defined by considering lines parallel
to the one through $(x,0)$ and $(0,1)$, and passing through $(0,y)$; or see \cite{descartes}, p.1, paragraph labelled 'la multiplication'.
\footnote{One can also interpret the field from $(V,\colin)$ alone, if $V$ is an $F$-vector space of dimension at least two, by similar means. 
 But as we are already given addition on $V$, and in particular on $\Qq a_1$, as well as a unit $a_1$, for our purposes it suffices to interpret multiplication.}%

 Now $(\Nn,+,\cdot)$ is interpretable in $(\Qq,+,\cdot)$ by Julia Robinson's theorem, using the theory of quadratic forms over $\Qq$.

 For the final nicety of avoiding the use of algebraic parameters, we will also need to know that $\Qq$ is interpretable in $(V,+,\colin)$ where $V$ is
 a $\Qq$-vector space of dimension $\geq 2$. 
 For this we take a field element to be a formal ration $u/v$ where $v \neq 0$ and $u \in \Qq v$ (we assume $\Qq v$ is definable uniformly in $v$, i.e. colinearity is definable). 
 If $v,v'$ are nonzero and noncolinear, we identify $u/v$ with $u'/v'$ if the lines $(u,u')$ and $(v,v')$ are parallel. Addition and multiplication of such fractions are definable in a similar way to the above.
 \vspace{1mm} \newline
\noindent \emph{Restricting to $\kk_0(t)$-parameters.} Let $G=\Gal(L/k_0(t))$. Let $Y=Gb_1 \union G b_2 $ be the (finite) set of Galois conjugates of $b_1,b_2$. Since $C_1$ is Galois over the $t$-line, all the elements of $Y$ are in $C_1(L)$. The vector subspace of $H$
generated by the images of the elements of $Y$ is still $\kk_0(t)$-definable and finite-dimensional, but at least two-dimensional; and is interpretable in $(\cti,t)$ along with colinearity. So again the coefficient field is interpretable.
\end{proof}

\begin{prop}\label{withoutt} Without the parameter $t$, the complexity of $\Th(\cti)$ is at least that of full arithmetic, i.e. $0^{(\omega)}$; and thus equal to 
$0^{(\omega)}$ if $\kk$ is countable. 
 \end{prop}

\prf It suffices to find a nonempty definable set $\Delta$ and a ring $\Psi(a) = (\Psi(a), +_a,\cdot_a)$ interpreted uniformly in $a \in \Delta$,
such that for any $a \in \Delta$, $\Psi(a)(\cti) \cong (\Zz,+,\cdot)$. For in this case a sentence $\theta$ holds true in $\Zz$ iff $\cti \models (\forall a \in \Delta)(\theta^{\Psi(a)})$,
where $\theta^U$ relativizes the quantifiers of $\theta$ to $U$.

Let $\Delta_1$ be the set of nonzero $a \in H$ represented by an irreducible divisor on $C$ for some genus one curve $C$. For $a \in \Delta_1$, 
we saw in \lemref{1d} that $\Phi(a) = \Qq a$, and also saw how to interpret a field $\Psi(a) $ and an action of $\Psi(a)$ on $\Qq a$, so that $\Psi(a)$ becomes isomorphic
to $\End_{\Qq} \Phi(a)$; this via formulas $\Psi,+,\cdot$ that are uniform in $a$.

Let $\Delta$ be the set of all elements $a \in H \m (0)$ such that the formulas $\Psi(a),+_a,\cdot_a$ define a field acting on $\Phi(a)$. Then $\Delta_1 \subset \Delta$, hence $\Delta$ is nonempty. Now by \remref{1dr}, for all $a \in H$, $\Phi(a) $ is a finite-dimensional $\Qq$-space. Hence for all $a \in \Delta$, $\Psi(a)$ is a number field.
By \cite[Theorem 3, pg 201]{rumely} there exists a formula uniformly defining $\Zz$ within any number field, finishing the proof. \eprf

\ssec{Questions}

 The call on Rumely's theorem in \propref{withoutt} can be avoided if we find a nonempty, zero-definable subset of $H$ where $\Phi(a)= \Qq a$. It seems likely that
this is the case, possibly for all $a$, or at least for $a$ such that for all $b \in \Bb_f$, the intersection $H(b) \meet \Phi(a)$ is either $0$ or $\Phi(a)$. 
It would be interesting to determine this.

Consider the structure $\mS$ described above; view the support relation on ${{H}} \times B$ as a basic relation
(asserting of $(h,b)$ that there exists $d \in D$ with $\rho(d)=h$ and supported on $b$.) Is the universal theory of $\mS$ decidable? 
Does it have a model companion? Is there an analogous decidable theory, where the Picard groups are taken with real coefficients? 

 The formula defining colinearity is {\em universal}: $x \in \Qq c$ iff for all $b \in B$, if $c$ is supported on $b$ then so is $x$.
 
Can $\Qq$ also be interpreted by means of universal formulas?

 Can we also interpret $\Qq$ within ${{\Bb_f}}$? 
 
 Is the theory of the ring $\cti$ recursively axiomatizable relative to $\Bb,{H}$?

 The theory of $k[t]^{int}$ depends at most on the characteristic of $k$, and on the
 transcendence degree over the prime field. Does it in fact depend on the latter? 
(In case the constant field $k$ is definable, the answer would be yes, and one has a strong form of undecidability
 by interpreting all finite subsets of $k$. 
 See \cite{raphael-robinson} for similar statements. 
 Also, \cite{poonen} may be relevant here.)

 Are the various fields $\Psi(a)$, $a \in \Delta_1$ definably isomorphic? (If so one can interpret $\Zz$ without parameters; but the question seems geometrically interesting on its own.)

\end{section}

\end{document}